\documentstyle[12pt,amssymb]{article}

   \makeatletter
   \ifnum\@ptsize=0  \fi
   \ifnum\@ptsize=2  \fi
   \if@twoside\fi
   \def\maketitle{\begin{center}\let\thanks=\footnote
      {\large\bf\@title}\par\bigskip\bigskip{\sc\@author}\par\bigskip{\rm\@date}
      \end{center}\bigskip\thispagestyle{empty}}

   \def\@begintheorem#1#2{\trivlist\item[\hskip\labelsep{\bf#1 #2.}]\it}
   \def\@opargbegintheorem#1#2#3{%
      \trivlist\item[\hskip\labelsep{\bf#1 #2 {\rm(#3)}.}]\it}

   \let\o@item=\@item
   \def\@item[#1]{\o@item[\rm #1]}

   \def\@biblabel#1{#1.}
   \let\othebibliography=\thebibliography
   \def\thebibliography#1{\small
      \def\@listi{\topsep=0cm\parsep=0cm\itemsep=0cm}\othebibliography{#1}}

   \newenvironment{commdiag}{\arraycolsep=0.1em
      \begin{array}{cccccccccccccc}%
   }{\end{array}}

   \newlength{\rightarrowlength} \rightarrowlength=2.5em
   \newlength{\downarrowlength}  \downarrowlength=5ex
   \def\Rarrow{\mathop{\makebox[\rightarrowlength]{\rightarrowfill}}\limits}
   \def\Darrow{\left\downarrow\parbox{0cm}{\rule{0cm}{\downarrowlength}}\right.}
   \def\Rhook{\lhook\mkern-9mu\Rarrow}
   \def\Dhook{%
      \left\downarrow\parbox{0cm}{\rule{0cm}{\downarrowlength}}\right.%
      \kern-1.46ex\raise0.50\downarrowlength\hbox{$\scriptscriptstyle\cap$}%
   }
   \def\Dequal{\left\|\parbox{0cm}{\rule{0cm}{\downarrowlength}}\right.}
   \def\rlabel#1{\rlap{$\scriptstyle#1$}}
   \def\llabel#1{\llap{$\scriptstyle#1$}}

   \newtheorem{satz}{Satz}[section]
   \makeatletter
   \@addtoreset{equation}{satz}
   
   \newtheorem{theorem}[satz]{Theorem}
   
   \newtheorem{remark}[satz]{Remark}
   \newtheorem{lemma}[satz]{Lemma}
   \newtheorem{corollary}[satz]{Corollary}

   \newenvironment{varthm*}[1]{\trivlist\item[]\it{\bf #1.}}{\endtrivlist}
   \def\startproof{\addvspace{\bigskipamount}\noindent}
   \def\proof{\startproof{\it Proof. }}
   \def\proofof#1{\startproof{\it Proof of #1.}}
   \def\qedsymbol{\frame{\rule[0pt]{0pt}{8pt}\rule[0pt]{8pt}{0pt}}}
   \def\qed{\nopagebreak\hspace*{\fill}\qedsymbol\par\addvspace{\bigskipamount}}

   \def\to{\longrightarrow}
   \def\mapsto{\mapstochar\longrightarrow}
   \def\epsilon{\varepsilon}
   \def\tilde{\widetilde}
   \def\hat{\widehat}
   \def\bar{\overline}
   \def\({\left(}
   \def\){\right)}
   \def\O{{\cal O}}

   \def\subjclass#1{{\def\thefootnote{}
      \footnotetext{1991 {\it Mathematics Subject Classification:} #1.}}}

   \newenvironment{items}{\list{$\bullet$}{
      \parsep=0cm\itemsep=0cm\topsep=0cm\partopsep=\medskipamount
      \def\makelabel##1{\hss\llap{##1}}}
   }{\endlist}

   \def\bbP{{\Bbb P}}

   \def\bbR{{\Bbb R}}
   \def\bbC{{\Bbb C}}
   \def\bbZ{{\Bbb Z}}
   \def\bbone{{\mathchoice {\rm 1\mskip-4mu l} {\rm 1\mskip-4mu l}
      {\rm 1\mskip-4.5mu l} {\rm 1\mskip-5mu l}}}
   \def\be{\begin{eqnarray*}}
   \def\ee{\end{eqnarray*}}
   
   \def\Bigwith{\ \Big\vert\ }
   \def\liste#1#2#3{\mbox{$#1_{#2},\dots,#1_{#3}$}}
   \def\eqnref#1{(\ref{#1})}
   \def\ds{\displaystyle}
   \def\tensor{\otimes}
   
   \def\with{\mid}
   \def\isom{\cong}
   \def\inverse{^{-1}}
   \def\Cal#1{{\cal #1}}
   \def\vect#1#2{\left(\begin{array}{c}#1\\ #2\end{array}\right)}
   \def\operatorname#1{\mathop{\rm #1}\nolimits}
   \def\diag{\operatorname{diag}}

   \def\Im{\operatorname{Im}}
   \def\NS{\operatorname{NS}}

   \def\rank{\operatorname{rank}}
   \def\eqdef{=_{\operatorname{def}}}
   \def\mult{\operatorname{mult}}

   \def\and{\quad\mbox{ and }\quad}
   
   \def\eps{\epsilon}
   \def\Sp{\operatorname{Sp}}
   \def\Bl{\operatorname{Bl}}
   \def\vol{\operatorname{vol}}
   \def\transp#1{{}^t#1}
   \def\matr#1#2#3#4{\left(\begin{array}{cc}#1 & #2 \\ #3 & #4
      \end{array}\right)}
   \def\symplD{\matr 0D{-D}0}
   \def\stackit#1#2{\stackrel{\scriptstyle #1}{#2}}
   \def\siegel{{\frak{H}}}
   \def\Nm{N\hskip-0.3ex m}

   \let\oldint=\int   \def\int{\oldint\limits}

\begin{document}

\title{Seshadri constants and periods of polarized abelian varieties}
\author{Thomas Bauer}
\date{}
\maketitle


\setcounter{section}{-1}
\section{Introduction}

   The purpose of this paper is to study the Seshadri constants of
   abelian varieties.  Consider 
   a polarized abelian variety $(A,L)$ 
   of dimension $g$ over the field
   of complex numbers.  
   One can associate to $(A,L)$ a real number $\eps(A,L)$, its {\em
   Seshadri constant}, which in effect measures how much of the
   positivity of $L$ can be concentrated at any given point of $A$.
   The number $\eps(A,L)$ can be defined as the 
   rate of growth in $k$ of the number
   of jets that one can specify in the linear series $|\O_A(kL)|$.
   Alternatively, one considers the blow-up
   $
      f:\tilde X=\Bl_x(X)\to X
   $ 
   of $X$ at a point $x$ with exceptional divisor 
   $E\subset\tilde X$ over $x$, and defines
   $$
      \eps(A,L)\eqdef
	 \sup\{\ \eps\in\bbR\with f^*L-\eps E\mbox{ is nef }\} \ .
   $$
   (Since $A$ is homogeneous, this is independent of $x$.)
   There has been recent interest in finding bounds on the Seshadri
   constants of abelian varieties and on smooth projective varieties
   in general (see \cite{EinKueLaz95} and \cite{Laz96}).  
   For the case of abelian varieties
   one has the elementary bounds
   $$
      1 \le\eps(A,L)\le\sqrt[g]{L^g} \ ,
   $$
   where by a result of Nakamaye \cite{Nak96} the
   lower bound is taken on only by abelian varieties which are
   polarized products of an elliptic curve and an abelian
   variety of dimension $g-1$.  
   
   Write now, as usual, $A$ as
   the quotient $A=V/\Lambda$ of its universal covering $V$ and a
   lattice $\Lambda\subset V$.  Viewing the first Chern class of $L$
   as a positive definite Hermitian form on $V$, its real part
   is a positive definite inner product $b_L$ on $V$, where we
   consider $V$ as a real vector space of dimension $2g$.  We
   define the {\em minimal period length} of $(A,L)$ to be
   real number
   $$
      m(A,L)\eqdef\min_{\stackit{\lambda\in\Lambda}{\lambda\ne 0}}
      b_L(\lambda,\lambda) \ .
   $$
   So $m(A,L)$ is the (square of the) length of the shortest
   non-zero period of $A$, where the length is taken with
   respect to the euclidian metric defined by $b_L$.  
   When $L$ is a principal polarization, this invariant has been
   studied by Buser and Sarnak in \cite{BusSar94}, who
   use an average argument
   familiar from the geometry of numbers to get a bound
   on the maximal value of
   $m(A,L)$.
   Lazarsfeld has recently established in \cite{Laz96} a
   surprising connection between minimal period lengths and
   Seshadri constants.  Using symplectic blowing up in the spirit of
   \cite{McDPol94} he shows that the Seshadri constant of
   $(A,L)$ is bounded below in terms of $m(A,L)$:
   $$
      \eps(A,L)\ge \frac{\pi}4 m(A,L) \ . \eqno(L)
   $$
   By generalizing the result of Buser and Sarnak
   we obtain a lower bound on 
   $m(A,L)$ in terms of the type of the polarization, 
   which then combined with $(L)$ leads to:

\def\rmitem[#1]{{\rm#1} }
\bgroup\def\thesatz{\arabic{satz}}
\begin{theorem}
   \rmitem[(a)]
   The maximal value of $m(A,L)$ as $(A,L)$ varies over the moduli
   space $\Cal A_D$ of polarized abelian varieties of fixed type
   $(\liste d1g)$ is bounded below by
   $$
      \max_{(A,L)\in\Cal A_D} m(A,L) \ge \frac1{\pi}\sqrt[g]{2L^g} \
      .
   $$

   \rmitem[(b)]
   For the very general polarized abelian variety $(A,L)$ of
   fixed type $(\liste d1g)$ one has the inequality
   $$
      \eps(A,L)\ge\frac14\sqrt[g]{2L^g}=
      \frac14\(2g!\prod_{i=1}^g d_i \)^{\frac1g} \ .
   $$
\end{theorem}
   Note that for $g\gg 0$ the upper bound given in the theorem
   differs from the theoretical upper bound $\sqrt[g]{L^g}$ only by a
   factor of approximately $4$.  Also note that in \cite{Laz96} the
   inequality $(L)$ 
   is stated for principally polarized abelian varieties.
   The proof given there, however, extends immediately to abelian
   varieties with polarizations of arbitrary type.

   It is well-known that bounds on Seshadri constants have
   implications for adjoint linear series.  In our situation this
   applies to the question whether an ample line bundle of some type
   $(\liste d1g)$ is very ample.  One knows by the classical theorem
   of Lefschetz that $L$ is very ample whenever $d_1\ge 3$.  Further,
   in case $d_1=2$ a result of Ohbuchi states that $L$ is very ample
   if and only if the linear series $|\O_A(\frac12L)|$ has no fixed
   divisor.  The remaining case of primitive line bundles,
   i.e.\ those with $d_1=1$, seems however to be very hard to deal
   with.  Debarre, Hulek and Spandaw consider in \cite{DebHulSpa94}
   polarizations of type $(1,\dots,1,d)$, i.e.\ pullbacks of
   principal polarizations under cyclic isogenies, and show that for
   a generic $(A,L)$ of this type, $L$ is very ample as soon as
   $d>2^g$.  Theorem 1 implies a criterion of a similar flavor for
   polarizations of arbitrary type, although the actual number that
   one gets in the special case of type $(1,\dots,1,d)$ is worse:

\begin{corollary}
   Let $(A,L)$ be a generic polarized abelian variety of type
   $(\liste d1g)$.  If
   $$
      \prod_{i=1}^g d_i
      \ge\frac{(8g)^g}{2g!}\approx\frac12(8e)^g \ ,
   $$
   then $L$ is very ample.
\end{corollary}

   Roughly speaking, Theorem 1 says that the Seshadri
   constant of a very general abelian variety $(A,L)$ is quite close
   to the theoretical upper bound $\sqrt[g]{L^g}$.  In the
   other direction, one is lead to ask under which geometrical
   circumstances $\eps(A,L)$ can become small -- apart from
   the trivial situation when $A$ contains an elliptic curve
   of small degree.  Lazarsfeld \cite{Laz96} has shown that
   for the Jacobian $(JC,\Theta)$ of a compact Riemann surface
   $C$ of genus $g\ge 2$ one has $\eps(JC,\Theta)\le\sqrt{g}$.
   Now the principally polarized abelian varieties which may
   be considered as being closest to Jacobians are Prym
   varieties of \'etale double coverings.  Our second result then
   shows that this intuition is indeed reflected by the fact that Prym
   varieties have small Seshadri constants:

\begin{theorem}
   Let $(P,\Xi)$ be the Prym variety of an \'etale double covering
   $\tilde C\to C$ of a compact Riemann surface $C$ of genus
   $g\ge 3$.
\begin{items}
   \item[(a)]
      One has
      $$
	 \eps(P,\Xi)\le\sqrt{2(g-2)}=\sqrt{2(\dim(P)-1)} \ .
      $$
   \item[(b)]
      If $C$ admits a map $C\to\bbP^1$ of degree $d$, and if
      $\tilde C$ is not hyperelliptic, then
      $$
         \eps(P,\Xi)\le\frac{2(d-1)(g-1)}{d+g-1}
	 =\frac{2(d-1)\dim(P)}{d+\dim(P)} \ .
      $$
\end{items}
\end{theorem}

   As for the assumptions in (b) note that for hyperelliptic
   $\tilde C$, $(P,\Xi)$ is a hyperelliptic Jacobian, and then
   one has $\eps(P,\Xi)\le 2\dim(P)/(1+\dim(P))$ by
   \cite{Laz96}, which is weaker than the inequality for the
   non-hyperelliptic case in (b).  Already the case $\dim(P)=2$
   shows however that this cannot be improved in general.

   The main result of the paper \cite{BusSar94} by Buser and Sarnak
   states that Jacobians have periods of unusually small length.  It
   is a consequence of Theorem 3 that a similar statement also holds
   (with larger numbers) for Prym varieties.  In fact, combining
   Theorem 3 with $(L)$ we obtain:

\begin{corollary}
   In the situation of cases (a) and (b) of Theorem 3 one has
   the following bounds on the minimal period length
   $m(P,\Xi)$:
\begin{items}
   \item[(a)]
      $\displaystyle m(P,\Xi)\le\frac4{\pi}\sqrt{2(g-2)}$,

   \medskip
   \item[(b)]
      $\displaystyle m(P,\Xi)\le\frac{8(d-1)(g-1)}{\pi(d+g-1)}
      \le\frac{8(d-1)}{\pi}$.
\end{items}
\end{corollary}

   It would be interesting to know if one can get stronger
   inequalities by the methods of \cite{BusSar94}.

   Finally, in an appendix (joint with T.\ Szemberg) we show how one
   can obtain more refined results on Seshadri constants for the case of
   abelian surfaces.  One knows by work of Steffens \cite{Ste} that
   for an abelian surface $(A,L)$ 
   of type $(1,d)$ the Seshadri constant is
   maximal, i.e.\ equal to $\sqrt{2d}$, if $2d$ is a square and 
   $\rank\NS(A)=1$.  The most surprising result here is that by
   contrast if $2d$ is not a square, then $\eps(A,L)$ is always
   sub-maximal:

\begin{theorem}[with Szemberg]
   Let $A$ be an abelian surface and let $L$ be an ample line bundle
   of type $(1,d)$, $d\ge 1$.  
   If $\sqrt{2d}$ is irrational, then
   $$
      \eps(L) \le \frac{2d}{\sqrt{1/k_0^2+2d}} \ ,
   $$
   where $(\ell_0,k_0)$ is the primitive solution of the diophantine
   equation
   $\ell^2 - 2dk^2 = 1$ (Pell's equation).
   In particular $\eps(L)$ is sub-maximal, i.e.\
   $\eps(L) < \sqrt{2d}$.

   If $2d+1$ is a square, then the inequality above is sharp.  In
   fact, in this case the upper bound is taken on whenever
   $\NS(A)\isom\bbZ$.
\end{theorem}

   As a consequence, one obtains:

\begin{corollary}
   The Seshadri constant of an ample line bundle on an abelian
   surface is rational.  
\end{corollary}
\egroup

   It is not known if Seshadri constants are always rational numbers,
   not even for the case of abelian varieties or for smooth
   surfaces.  Also, 
   the sub-maximality statement in Theorem 5 suggests the
   possibility that there is additional structure to these invariants
   which is not fully understood yet.


\begin{varthm*}{\it Acknowledgements}
\rm
   This research was done during the author's stay at the
   University of California, Los Angeles, which was supported
   by DFG grant Ba 1559/2-1. 
   It is my pleasure to thank R.\ Lazarsfeld for many helpful
   discussions and UCLA for its hospitality.
\end{varthm*}

\begin{varthm*}{\it Notation and Conventions}
\rm
   We work throughout over the field $\bbC$ of complex numbers.

   Numerical equivalence of divisors or line bundles will be
   denoted by $\equiv$.
\end{varthm*}


\section{Period lengths of abelian varieties}
   \label{minimal period lengths}

   The purpose of this section is to prove Theorem 1 from the
   introduction.  We start with some remarks on polarized abelian
   varieties.
   So let $\liste d1g$ be positive integers such that
   $d_i|d_{i+1}$ for $1\le i<n$, let $D$ be the diagonal
   matrix $D=\diag(\liste d1g)$ and denote as usual by $\Cal
   A_D$ the moduli space of polarized abelian varieties of
   type $(\liste d1g)$.  Recall that it can be realized as a
   quotient 
   $$
      \Cal A_D=\siegel_g/\Sp_{2g}^D(\bbZ)
   $$
   of the Siegel upper half space 
   $ \siegel_g=\{Z\in M_g(\bbC)\with \transp Z=Z\mbox{ and }\Im Z>0\} $
   where the symplectic group
   $$
      \Sp_{2g}^D(\bbZ)=\left\{R\in M_{2g}(\bbZ)\Bigwith 
      R\symplD\transp R=\symplD\right\} \ ,
   $$
   acts on $\siegel_g$ by
   $$
      Z\mapsto R\cdot Z=(aZ+bD)(D\inverse cZ+D\inverse
      dD)\inverse\qquad\mbox{ for } 
      R=\matr abcd\in\Sp_{2g}^D({\bbZ})
   $$
   (cf.\ \cite[Sect.\ 8.2]{LB}).
   Following the approach of Buser and Sarnak \cite{BusSar94} we now
   show:

\begin{theorem}\label{thm period lengths}
   \rmitem[(a)]
   One has
   $$
      \max_{(A,L)\in\Cal A_D} m(A,L) \ge \frac1{\pi}\sqrt[g]{2L^g} \
      .
   $$

   \rmitem[(b)]
   There is a countable union $\Cal B\subset\Cal A_D$ of
   proper closed subvarieties of $\Cal A_D$ such that for all
   $(A,L)\in\Cal A_D-\Cal B$ one has the inequality
   $$
      \eps(A,L)\ge\frac14\sqrt[g]{2L^g}=
      \frac14\(2g!\prod_{i=1}^g d_i \)^{\frac1g} \ .
   $$
\end{theorem}

\proof
   Let $Z$ be an element of $\siegel_g$, i.e.\ $Z=X+iY$ with
   real-valued symmetric matrices $X$ and $Y$ such that 
   $Y$ is positive definite.  Recall that modulo the action of
   $\Sp_{2g}^D(\bbZ)$ on $\siegel_g$ the matrix $Z$ corresponds
   to the isomorphism class of the 
   polarized abelian variety $(A_Z,L_Z)$ whose lattice
   is
   $$
      \Lambda_Z=(Z,D)\bbZ^{2g}\subset\bbC^g
   $$
   and whose Hermitian form $H_Z$ is given by the matrix
   $Y\inverse$ with respect to the standard basis of $\bbC^g$.
   Consider the quadratic form
   \be
      q_Z:\Lambda_Z&\to&\bbR \\
	  \lambda&\mapsto& H_Z(\lambda,\lambda) \ .
   \ee
   The columns of the matrix $(Z,D)$ form a symplectic basis
   for $\Lambda_Z$, i.e.\ a basis with respect to which the
   first Chern class of $L_Z$, viewed as an alternating form
   on $\Lambda_Z$, is given by the matrix $\symplD$.  Now for
   $m,n\in\bbZ^g$ we have
   \be
      q_Z(Zn+Dm) &=& 
      \transp{(\bar{Zn+Dm})} Y\inverse (Zn+Dm) \\
      &=& (\transp m, \transp n)\matr D0XD
      \matr{Y\inverse}00{D\inverse Y D\inverse}\matr
      DX0D\vect mn
   \ee
   so that $q_Z$ is
   given with respect to the symplectic basis
   by the matrix $Q_Z=\transp{P_Z}P_Z$, where
   \begin{equation}\label{def P} 
      P_Z=P_{X,Y}=\matr{\sqrt{Y\inverse}}00{\sqrt YD\inverse}
      \matr DX0D \in M_{2g}(\bbR) \ .
   \end{equation} 

   Fix now a real number $R>0$ and denote for a polarized
   abelian variety $(A,L)$ by $n_R(X,L)$ the number of non-zero
   periods in the closed ball $\bar{B_{R^2}(0)}\subset(V,b_L)$, i.e.
   $$
      n_R(A,L)= \#\{\lambda\in\Lambda-\{0\}\with
      b_L(\lambda,\lambda)\le R^2\} \ .
   $$
   In view of what we found above, one has
   $$
      n_R(A_Z,L_Z)=\sum_{\stackit{\ell\in\bbZ^{2g}}{\ell\ne 0}}
      \chi_{R^2}(\transp{\ell}Q_Z\ell) \ ,
   $$
   where $\chi_{R^2}$ is the characteristic function of the
   interval $[0,R^2]$.  We will now consider in particular
   period matrices of the form $Z=X+i\frac1{y^2}\bbone$ for
   $y>0$, where $\bbone$ denotes the identity matrix.  
   The idea is to 
   study the average of $n_R(A_Z,L_Z)$
   when $y$ is fixed and $X$ varies over a suitable compact
   set.  Specifically, let 
   $$
      V\subset\{X\in M_g(\bbR)\with\transp X=X\}
   $$
   be the compact subset consisting of the matrices whose
   entries are bounded by the exponent of $L_Z$, i.e.\
   \begin{equation}\label{def V} 
      V=\left\{X\Bigwith 0\le X_{ij}\le d_g\mbox{ for }1\le
      i,j\le g \right\} \ ,
   \end{equation} 
   and consider the average
   \be 
      I(y)&\eqdef&
	 \frac1{\vol(V)}
	 \int_V\sum_{\stackit{\ell\in\bbZ^{2g}}{\ell\ne
	 0}} f\(P_{X,\frac1{y^2}\bbone}\cdot \ell\) dX \\
	 &=& \frac1{\vol(V)}\int_V n_R(A_Z,L_Z) dX \ ,
   \ee
   where $f:\bbR^{2g}\to\bbR$ is the function
   $f(x)=\chi_{R^2}(\transp x\cdot x)$.
   It follows from Lemma \ref{calc} below that
   $$
      \lim_{y\to\infty} 
      I(y)=\frac{R^{2g}\cdot\sigma_{2g}}{\prod_{i=1}^g d_i} \ ,
   $$
   so that we will have $\lim\limits_{y\to\infty} I(y)<2$, if 
   we choose $R$ such that
   \begin{equation}\label{inequality R}
      R^2<\frac 1{\pi}\sqrt[g]{2g!\prod d_i} \ .
   \end{equation}
   But then there exists a real number $y>0$ and a
   symmetric matrix $X\in M_{2g}(\bbR)$ such that
   $$
      \sum_{\stackit{\ell\in\bbZ^{2g}}{\ell\ne 0}} 
	 f\(P_{X,\frac 1{y^2}\bbone}\cdot \ell\)
	 = n_R(A_Z,L_Z) < 2 \ .
   $$
   Since $n_R(A_Z,L_Z)$ is in any event an even non-negative
   integer, we must then have $n_R(A_Z,L_Z)=0$.  But this just
   means that for the polarized abelian variety $(A_Z,L_Z)$ 
   corresponding to $Z$ one has
   $$
      m(A_Z,L_Z) > R^2 \ ,
   $$
   and, using \eqnref{inequality R} this implies the asserted lower
   bound on the maximum of $m(A,L)$. 
   Assertion (b) follows from (a) and Lazarsfeld's inequality
   $(L)$.
\qed

\begin{lemma}\label{calc}
   Let $f:\bbR^{2g}\to\bbR$ be an integrable function of
   compact support.  Consider the function $I_f:\bbR^+\to\bbR$
   which is defined by
   $$
      I_f(y)=\frac 1{\vol(V)}
      \int_V\sum_{\stackit{\ell\in\bbZ^{2g}}{\ell\ne 0}}
      f\(P_{X,\frac1{y^2}\bbone}\cdot \ell\) dX 
   $$
   where $P$ and $V$ are as in \eqnref{def P} and \eqnref{def V}
   respectively.  Then
   $$
      \lim_{y\to\infty} I_f(y)=\frac1{\det
      D}\int_{\bbR^{2g}}f(x) dx \ .
   $$
\end{lemma}

   The proof is completely elementary but somewhat tricky.
   Here the  dependence on the type of the polarization comes in
   crucially.

\proofof{Lemma \ref{calc}}
   The integral $\vol(V)\cdot I_f(y)$ can be written as
   $$
      \int_V\sum_{\stackit{m,n\in\bbZ^g}{(m,n)\ne(0,0)}}
      f\left(\begin{array}{c}
	 y(d_1m_1+\sum_{i=1}^g X_{1i} n_i) \\
	 \vdots \\
	 y(d_gm_g+\sum_{i=1}^g X_{gi} n_i) \\
	 y\inverse n_1 \\
	 \vdots \\
	 y\inverse n_g
      \end{array}\right)
      dX \ .
   $$
   We integrate under the sum and consider first the terms
   with $n\ne 0$.  The contribution of such a term, if say
   $n_k\ne 0$, is
   $$
      c_n=\int_0^{d_g}\dots\int_0^{d_g}\sum_{m_g}\dots\sum_{m_2}
      F_{1,k}(m,n,X) dX' \ ,
   $$
   where we set $dX'=\prod_{\stackit{i\le j}{(i,j)\ne(1,k)}}
   dX_{ij}$ and
   \be
      F_{1,k}(m,n,X)&=&\sum_{m_1}\int_0^{d_g}
      f\vect{yd_1n_k(\frac{m_1}{n_k}+\frac{X_{1k}}{d_1}+\lambda_k)}{\vdots}
      dX_{1k} \\
      &=&\sum_{m_1=0}^{n_k-1}\sum_{j=-\infty}^{\infty}\int_0^{d_g/d_1}
      d_1
      f\vect{yd_1n_k(\frac{m_1}{n_k}+j+T_{1k}+\lambda_k)}{\vdots}
      dT_{1k} \\
   \ee
   where $\lambda_k$ is independent of $T_{1k}$ and $m_1$.  We
   therefore obtain
   \be
      F_{1,k}(m,n,X) 
      &=&n_k\frac{d_g}{d_1}\int_{-\infty}^{\infty} d_1 
      f\vect{yd_1n_kT_{1k}}{\vdots} dT_{1k} \\
      &=&y\inverse\frac{d_g}{d_1}\int_{-\infty}^{\infty}
      f\vect{t_1}{\vdots} dt_1 \ .
   \ee
   Continuing in the same manner with $m_2,X_{2k}$ up to
   $m_g,X_{2g}$ we find
   \be 
      c_n&=&\int_0^{d_g}\dots\int_0^{d_g}
      y\inverse\left(\prod_{i=1}^g\frac{d_g}{d_1}\right)
      \int_{\bbR^g} f(\liste t1g,\liste{y\inverse n}1g) dt
      \prod_{\stackit{i\le j}{i\ne k\ne j}} dX_{ij} \\
      &=& \frac{\vol(V)}{\det D} y\inverse\int_{\bbR^g} f(\liste
      t1g,\liste{y\inverse n}1g) dt
   \ee
   so that, taking into account that $f$ is of compact support,
   one gets
   $$
      \lim_{y\to\infty}\vol(V)\cdot I_f(y)
      =\lim_{y\to\infty}\sum_{n\in\bbZ^g} c_n + \sum_{m\ne 0}
      f\vect{yDm}{0}
      =\frac{\vol(V)}{\det D}\int_{\bbR^{2g}} f(x) dx
   $$
   which proves the lemma.
\qed

\begin{corollary}
   Let $(A,L)$ be a generic polarized abelian variety of type
   $(\liste d1g)$.  If
   $$
      \prod_{i=1}^g d_i \ge \frac{(8g)^g}{2g!} \approx
      \frac12 (8e)^g \ ,
   $$
   then $L$ is very ample.
\end{corollary}

   In fact, the bound on $\prod_{i=1}^g d_i$ guarantees by Theorem
   \ref{thm period lengths} that for the very general polarized
   abelian variety $(A,L)$ of the given type one has $\eps(A,L)\ge
   2g$.  This implies by a standard application of Kawamata-Viehweg
   vanishing (cf.\ \cite[Sect.\ 4]{EinKueLaz95} and \cite[Proposition
   6.8]{Dem92}) that $L$ is very ample.  Note that since very
   ampleness is an open condition on $\Cal A_D$, the corollary holds
   for {\em generic} $(A,L)$, even if we have the lower bound on
   $\eps(A,L)$ only for {\em very general} $(A,L)$.


\section{Seshadri constants of Prym varieties}

   Let $f:\tilde C\to C$ be an \'etale double cover of a
   compact Riemann surface $C$ of genus $g$. 
   Identify as usual the Jacobians $JC$ and $J\tilde C$ with
   their respective dual abelian varieties and
   consider the
   pullback map $f^*:JC\to J\tilde C$.
   The {\em Prym variety} $P$ of the given double cover
   is the complementary abelian subvariety of
   the image of $f^*$
   in $J\tilde C$ (see \cite[Chap.\ 12]{LB} and \cite{Mum74}).  
   The canonical principal polarization
   $\O_{J\tilde C}(\tilde\Theta)$ on $J\tilde C$ restricts to
   twice a principal polarization $\O_P(\Xi)$ on the
   $(g-1)$-dimensional abelian variety $P$.  We prove in this
   section the following bounds on the Seshadri constant of
   $(P,\Xi)$:

\begin{theorem}\label{Prym bounds}
   Assume $g\ge 3$.  Then:
\begin{items}
   \item[(a)]
      One has
      $$
	 \eps(P,\Xi)\le\sqrt{2(g-2)} \ .
      $$
   \item[(b)]
      If $C$ admits a map $C\to\bbP^1$ of degree $d$, and if
      $\tilde C$ is not hyperelliptic, then
      $$
	 \eps(P,\Xi)\le\frac{2(d-1)(g-1)}{d+g-1} \ .
      $$
\end{items}
\end{theorem}

   Corollary 4 in the introduction follows from the theorem
   and Lazarsfeld's result $(L)$.

\proofof{Theorem \ref{Prym bounds}}
   (a) Note first that we may assume $\dim(P)\ge 3$, since for
   $\dim(P)=2$ the inequality is clear.  Further, we may assume
   that $\tilde C$ is not hyperelliptic:  in fact, otherwise
   $(P,\Xi)$ is a Jacobian (see \cite[Corollary 12.5.7]{LB})
   and then one has 
   $\eps(P,\Xi)\le(\dim(P))^{\frac12}\le(2\dim(P)-2)^{\frac12}$.

   We will as usual identify $J\tilde C$ and $P$ with their
   respective dual abelian varieties via the isomorphisms
   defined by the principal polarizations $\O_{J\tilde
   C}(\tilde\Theta)$ and $\O_P(\Xi)$.  The dual map of the
   inclusion $\iota:P\hookrightarrow J\tilde C$ gives then a
   surjective morphism $\hat\iota:J\tilde C\to P$.  Note that
   the composition $\iota\circ\hat\iota$ is just the norm
   endomorphism $N_P$ of $P$ (cf.\ \cite[Sect.\ 12.2]{LB}).
   We will study the image $S\subset P$ of the composed map
   $$
      \psi:\tilde C\times\tilde C\stackrel{\tilde s}{\to}
      J\tilde C\stackrel{\hat\iota}{\to} P \ ,
   $$
   where $\tilde s$ is the subtraction map
   $(x,y)\mapsto\O_{\tilde C}(x-y)$.  We verify first that 
   $$
      \dim(S)=2 \ .
   $$
   In fact, we have
   $$
      S=\psi(\tilde C\times\tilde C)=N_P(\tilde C)-N_P(\tilde
      C) \ ,
   $$
   so if $S$ were a point, then $\tilde C$ would be contained
   in the kernel of $\hat\iota$ which is certainly impossible.
   And if $S$ were a curve, then one would have $N_P(\tilde
   C)-N_P(\tilde C)=N_P(\tilde C)$, so that $N_P(\tilde C)$
   would be an elliptic curve; but this would imply
   $N_P(\tilde\Theta)=N_P(\tilde C+\dots+\tilde C)=N_P(\tilde
   C)$, contradicting the fact that the fibres of $N_P$ are of
   dimension $g(\tilde C)-\dim(P)\le g(\tilde C)-3$.

   We claim next that
   \begin{equation}\label{deg formula}
      \deg_{\Xi}(S)=\Xi^2\cdot S=\frac8{\deg(\psi)}(g-1)(g-2)
      \ .
   \end{equation}
   For the proof of \eqnref{deg formula}
   recall first that 
   $2\tilde\Theta\equiv \Nm^*\Theta+\hat\iota^*\Xi$, 
   where $\O_{JC}(\Theta)$ is the
   canonical principal polarization on $JC$ and $\Nm:J\tilde
   C\to JC$ is the norm map associated with $f$ 
   (see \cite[Proposition 12.3.4]{LB}).  So we have
   $$
      \psi^*\Xi\equiv 2\tilde s^*\tilde\Theta-\tilde s^*
      \Nm^*\Theta \ .
   $$
   Let now $F_1,F_2\subset C\times C$ be fibres of the two
   projections and $\Delta\subset C\times C$ the diagonal.
   For the pullback $s^*\Theta$ of $\Theta$ under the
   subtraction map $s:C\times C\to JC$ one has
   $s^*\Theta\equiv(g-1)(F_1+F_2)+\Delta$ (cf.\ \cite[Theorem
   4.2]{Rai89}).  So, using the commutative diagram
   $$
   \begin{commdiag}
      \tilde C\times\tilde C & \Rarrow^{f\times f} & C\times C \\
      \Darrow\llabel{\tilde s} & & \Darrow\rlabel{s} \\
      J\tilde C & \Rarrow^{\Nm} & JC
   \end{commdiag}
   $$
   we find
   $$
      \tilde s^* \Nm^*\Theta\equiv(f\times f)^* s^*\Theta
      \equiv 2(g-1)(\tilde F_1+\tilde
      F_2)+\tilde\Delta+\Gamma_{\tau} \ ,
   $$
   where $\tilde F_1,\tilde F_2\subset\tilde C\times\tilde C$
   are fibres of the projections, $\tilde\Delta$ is the
   diagonal, and $\Gamma_{\tau}$ is the graph of the covering
   involution $\tau:\tilde C\to\tilde C$.
   We conclude that
   \begin{equation}\label{eqn pullback Xi}
      \psi^*\Xi\equiv(2g-2)(\tilde F_1+\tilde
      F_2)+\tilde\Delta-\Gamma_{\tau}
   \end{equation}
   and, using the fact that $\tau$ is fixed-point free, this
   implies with a calculation
   \be
      \deg_{\Xi}(S)=\Xi^2\cdot S 
      &=& \frac1{\deg(\psi)}\Xi^2\cdot\psi_*(\tilde
      C\times\tilde C)
      \\
      &=&\frac1{\deg(\psi)}(\psi^*\Xi)^2
      \\
      &=&\frac8{\deg(\psi)}(g-1)(g-2) \ ,
   \ee
   as claimed.

   The diagonal $\tilde\Delta$ is the scheme-theoretic
   inverse image of $0$ under $\psi$.  In fact, one has
   $N_P=1-\tilde\tau$, where $\tilde\tau$ is the involution
   on $J\tilde C$ induced by $\tau$, thus for $x,y\in\tilde C$
   $$
      \iota\circ\psi(x,y)=\iota\circ\hat\iota\circ\tilde
      s(x,y)
      =N_P\O_{\tilde C}(x-y)
      =\O_{\tilde C}(x-y-\tau(x)+\tau(y)) \ ,
   $$
   so that, since $\tilde C$ is not hyperelliptic,
   $\psi(x,y)=0$ implies $x=y$.  So we have
   $\psi\inverse(0)=\tilde\Delta$ set-theoretically and, using
   the fact that the Abel-Prym map 
   $$
      \tilde C\hookrightarrow J\tilde C\stackrel{\hat\iota}{\to}P
   $$ 
   is an embedding, one checks
   that this also holds scheme-theoretically.  
   Let now $\tilde P\to P$ be the blow-up of $P$ at $0$ with
   exceptional divisor $E$, and let $\tilde S=\Bl_0(S)$ be the
   proper transform of $S$.  One has a commutative diagram
   $$
   \begin{commdiag}
      \tilde\Delta & \Rarrow & \bbP C_0(S) & \Rhook &
      E\rlap{ $=\bbP T_0(P)$} \\
      \Dhook & & \Dhook & & \Dhook \\
      \tilde C\times\tilde C & \Rarrow^{\tilde\psi} & \tilde S
      & \Rhook & \tilde P\rlap{ $=\Bl_0(P)$} \\
      \Dequal & & \Darrow & & \Darrow \\
      \tilde C\times\tilde C & \Rarrow^{\psi} & S
      & \Rhook & P
   \end{commdiag}
   $$
   where $\bbP C_0(S)$ is the projective tangent cone of $S$
   at $0$.  So, using $\tilde\psi^* E=\tilde\Delta$, we obtain
   \begin{eqnarray}
      \mult_0(S) & = &
	 \int_E \O_E(1)\cdot[\bbP C_0(S)] \nonumber
	 \\ 
	 &=& \int_{\tilde P} \O_{\tilde P}(-E)\cdot [\bbP C_0(S)]
	 = - \int_{\tilde P} \O_{\tilde P}(E)^2 \cdot [\tilde S]
	 \nonumber
	 \\
	 & = & -\frac1{\deg(\psi)} \int_{\tilde C\times \tilde C}
	    \tilde\psi^*\O_{\tilde P}(E)^2 
         = -\frac1{\deg(\psi)}\cdot\tilde\Delta^2
	 \nonumber
	 \\
	 & = & \frac1{\deg(\psi)}(2g(\tilde C)-2)
	 = \frac4{\deg(\psi)}(g-1) \ .  \label{mult formula}
   \end{eqnarray}
   Now recall that by \cite[(6.7)]{Dem92} any singular subvariety of $P$
   leads to an upper bound on the Seshadri constant of
   $(P,\Xi)$.  Applying this to $S$ we find upon using
   \eqnref{deg formula} and \eqnref{mult formula}
   \be
      \eps(P,\Xi)\le\sqrt{\frac{\deg_{\Xi}(S)}{\mult_0(S)}}
      =\sqrt{2(g-2)} \ .
   \ee

   (b) Suppose now that there exists a map
   $\phi:C\to\bbP^1$ of degree $d$.  This
   implies that there is an effective divisor 
   $D\in|d(F_1+F_2)-\Delta|$, namely the
   closure of $\{(x,y)\with\phi(x)=\phi(y), x\ne y\}$.  It pulls back
   to an effective divisor
   $$
      (f\times f)^* D\in|2d(\tilde F_1+\tilde
      F_2)-\tilde\Delta-\Gamma_{\tau}| \ .
   $$
   Since by assumption $\tilde C$ is not hyperelliptic, we
   have again $\tilde\Delta=\psi\inverse(0)$
   scheme-theoretically, and therefore the $\bbR$-divisor
   $\psi^*\Xi-\eps(P,\Xi)\cdot\tilde\Delta$ is nef, so that
   $$
      \(\psi^*\Xi-\eps(P,\Xi)\cdot\tilde\Delta\)\cdot(f\times
      f)^* D \ge 0 \ .
   $$
   Upon using \eqnref{eqn pullback Xi} one obtains the asserted
   inequality for $\eps(P,\Xi)$.
\qed


\def\ds{\displaystyle}
\def\docases#1{\left\{\begin{array}{cl} #1 \end{array}\right.}
\def\m{{\frak m}}
\def\rank{\operatorname{rank}}
\def\rounddown#1{\left\lfloor#1\right\rfloor}


\section*{Appendix: Seshadri constants of abelian surfaces\\
   \rm Thomas Bauer and Tomasz Szemberg}
\def\thesection{A}
\setcounter{satz}{0}

   Our purpose here is to show how one can get more refined
   results on Seshadri constants for the case of abelian surfaces.  
   In particular it follows that, somewhat surprisingly,
   Seshadri constants on abelian surfaces are always rational.

   Consider an abelian surface $A$ and an ample line bundle
   $L$ on $A$.  Since $\eps(kL)=k\eps(L)$ for any integer
   $k>0$, we may assume that
   $L$ is primitive, i.e.\ of type $(1,d)$ for some integer
   $d\ge 1$.  Recall the elementary bounds
   $$
      1 \le\eps(L)\le\sqrt{2d} \ . \eqno (*)
   $$
   One knows moreover by \cite[Theorem 1.2]{Nak96} that
   $\eps(L) \ge \frac43$, 
   unless $A$ is a product of elliptic curves.  
   Further, if $\sqrt{2d}$ is rational and $\rank\NS(A)=1$, then
   by \cite{Ste} the Seshadri constant $\eps(L)$ is maximal,
   i.e.\ $\eps(L)=\sqrt{2d}$, which shows that the upper bound in $(*)$
   cannot be improved in general.  
   On the other hand, if $\sqrt{2d}$ is irrational, then our 
   result shows that one does have a better upper bound: 

\begin{theorem}\label{thm surfaces}
   Let $A$ be an abelian surface and let $L$ be an ample line bundle
   of type $(1,d)$, $d\ge 1$.  
   \begin{items}
   \item[(a)]
      If $\sqrt{2d}$ is irrational, then
      $$
         \eps(L) \le \frac{2d}{\sqrt{1/k_0^2+2d}} \ ,
      $$
      where $(\ell_0,k_0)$ is the primitive solution of the diophantine
      equation
      $\ell^2 - 2dk^2 = 1$ (Pell's equation).
      In particular $\eps(L)$ is sub-maximal, i.e.\
      $\eps(L) < \sqrt{2d}$.
   \item[(b)]
      One has the lower bound
      $$
         \eps(L) \ge \min\left\{\eps_0, \frac{\sqrt7}2 \sqrt d \right\} \ ,
      $$
      where $\eps_0$ is the minimal degree (with respect to
      $L$) of the elliptic curves in $X$.
   \item[(c)]
      If $2d+1$ is a square, then the inequality in (a) is
      sharp.  In fact, in this case the upper bound is taken on whenever
      $\NS(A)\isom\bbZ$.
   \end{items}
\end{theorem}

   At first sight the bound in (a) might appear non-constructive because it
   involves the primitive solution of Pell's equation.
   This solution, however, can be
   effectively computed via continued fractions (see \ref{examples}
   for the numerical values for polarizations of small degree).
   As for (b) note that
   it is inevitable that small values of $\eps(L)$ occur for non-simple
   abelian surfaces regardless how large the type of the polarization
   may be, since
   for any given integer $e\ge 1$ there are abelian surfaces $(A,L)$
   of arbitrarily
   high degree $L^2$
   containing an elliptic curve of degree $e$.
   We do not expect that the particular bound in (b) is
   optimal; on the other hand it is tempting, in
   view of (c), to wonder whether the bound
   in (a) might be sharp in general.

\proofof{Theorem \ref{thm surfaces}}
   (a) Since $\eps(L)$ is invariant under algebraic equivalence, we
   may assume $L$ to be symmetric.  Recall that for any $n\ge 1$ the
   space of sections of $\O_A(nL)$ admits a decomposition
   $$
      H^0(A,\O_A(nL))=H^0(A,\O_A(nL))^+ \oplus H^0(A,\O_A(nL))^-
   $$
   into the spaces of even and odd sections whose dimensions are given
   by the formula
   $$
      h^0(A,\O_A(nL))^{\pm}=
      2+\frac{n^2 d}2-\frac{n^{\mp}(\O_A(nL))}4 \ ,
   $$
   where $n^{\mp}(\O_A(nL))$ is the number of odd respectively even
   halfperiods of the line bundle $\O_A(nL)$ (cf.\
   \cite[Corollary 4.6.6]{LB} and \cite[Theorem 3.1]{Bau94}).  
   So for even multiples
   $n=2k$ of $L$ we have in particular $h^0(A,\O_A(2kL))^+=2+2dk^2$.
   On the other hand, since an even section vanishes in halfperiods
   to even orders, it is at most
   $$
      1 + 3 + \dots + (m-1) = \(\frac m2\)^2
   $$
   conditions on an even section to vanish at a fixed halfperiod $x$
   to an even order $m$.  This implies that 
   $$
      H^0(A,\O_A(2kL)\tensor\Cal I_x^m)^+ \ne 0
   $$
   provided that $m\le 2\sqrt{2dk^2+1}$.  Thus there
   exists an even divisor $D\in|2kL|^+$ with multiplicity
   $$
      \mult_x(D)\ge\rounddown{2\sqrt{2dk^2+1}}
   $$
   at $x$.  The crucial point is now to avoid the round-down in this
   expression.  To this end consider the diophantine equation
   $$
      \ell^2 - 2dk^2 = 1 \ ,
   $$
   a special case of {\em Pell's equation}.  Since 
   $\sqrt{2d}$ is by assumption irrational, Pell's equation has a
   primitive solution $(\ell_0,k_0)$.  We conclude that
   $$
      \eps(L) \le \frac{L\cdot D}{\mult_x(D)}=
         \frac{4dk_0}{2\sqrt{2dk_0^2+1}}=
         \frac{2d}{\sqrt{1/k_0^2+2d}}<
         \sqrt{2d} \ ,
   $$
   as claimed.

   (b)  Let $C\subset A$ be an irreducible curve of arithmetic
   genus $p_a(C)>1$, let $x\in C$ and $m=\mult_x(C)$.  
   We have to show that
   $$
      \frac{L\cdot C}{m} \ge \frac{\sqrt7}2 \sqrt d \qquad\mbox{
         for all } x \in A \ .
   $$
   First observe that for the geometric genus one has
   $$
      p_g(C) \ge 2 \ .
   $$
   In fact, suppose to the contrary that $p_g(C) \le 1$ and consider
   the normalization $N \to C$.  Abelian varieties do not
   contain any rational curves, so $p_g(C)=1$ and the composed map
   $$
      N \to C \hookrightarrow A
   $$
   is -- after possibly translating $C$ -- a homomorphism of abelian
   varieties, and hence an embedding, which is absurd.

   The adjunction formula and the inequality $p_a(C)-p_g(C) \ge
   {m \choose 2 }$ then yield
   $$
      m \le \sqrt{C^2 - \frac74} + \frac12 \ .
   $$
   Combining this bound on the multiplicities of irreducible curves
   with Hodge index gives
   $$
      \eps(L)\ge\inf_{L\cdot C} 
	 \left\{ \frac{L\cdot C}{\sqrt{\frac{(L\cdot C)^2}{2d} - \frac74}
         + \frac12} \right\}
   $$
   where the infimum is taken over the degrees $L\cdot C$ of the irreducible
   non-elliptic curves $C\subset A$.
   But the real-valued function
   $$
      f(t)=\frac t{\sqrt{\frac{t^2}{2d}-\frac74} + \frac12}
   $$
   takes on its minimum at $t_0=2\sqrt{7d}$ with minimal value
   $\frac{\sqrt7}2 \sqrt d$ at $t_0$.  This implies the assertion.

   (c) By assumption we have $2d+1=\ell^2$ for some integer
   $\ell\ge 1$.  Then $d$ is an even number, and after
   possibly replacing $L$ by another symmetric translate one
   has
   $$
      h^0(A,L)^+=\frac d2 + 1 \ .
   $$
   Since the line bundle $L$ is primitive, it has both even and odd
   halfperiods.  Thus we may choose an odd halfperiod $x$, so that
   the number of conditions on a section in $H^0(A,L)^+$ to vanish
   at $x$ to order $2p+1$ is $2+4+\dots+2p=p(p+1)$.  Therefore
   there exists a divisor $D\in|L|^+$ with
   $$
      \mult_x(D)\ge 2\rounddown{\frac12 \sqrt{2d+1}-\frac12} +
      1=\ell \ ,
   $$
   Suppose now that there is an irreducible curve
   $C\subset A$ with $L\cdot C/\mult_x(C) < 2d/\ell$.
   Then the assumption $\NS(A)\isom\bbZ$ implies that 
   $D$ is irreducible, so that 
   $C$ and $D$ intersect properly, and hence 
   $$
      L\cdot C=D\cdot C\ge\mult_x(D)\cdot\mult_x(C)
      > (L\cdot D)(L\cdot C)\frac{\ell^2}{(2d)^2} =
      (L\cdot C)\frac{2d+1}{2d} \ ,
   $$
   a contradiction.
\qed

   In all known examples Seshadri constants
   have turned out to be rational.
   While it is unclear if this is true in general,
   Theorem \ref{thm surfaces} implies:

\begin{corollary}
   The Seshadri constant of an ample line bundle on an abelian surface is
   rational.
\end{corollary}

   In fact, as shown in \cite{Ste}, this follows from the 
   submaximality statement in part (a) of Theorem 
   \ref{thm surfaces} by
   applying the Nakai-Moishezon criterion for $\bbR$-divisors
   \cite{CamPet90}.
   It would be interesting to have a more conceptual explanation for
   the sub-maximality statement in Theorem \ref{thm surfaces} and
   also for the rationality of Seshadri constants on abelian
   surfaces.

\begin{remark}\label{examples}\rm
   In order to convey some feeling for the numbers involved here, 
   we list for $1\le d\le 20$ the (truncated) numerical 
   values of the upper bound $\eps_{\rm
   upper}(d)$ given in
   part (a) of Theorem \ref{thm surfaces} along with the lower
   bound $\eps_{\rm lower}(d)$ from part (b) and the
   theoretical upper bound $\sqrt{2d}$.  Note that the 
   lower bound holds for {\em simple} $A$ only.
   Numbers in boldface indicate the cases where one knows the
   exact value of $\eps(A)$ when $\rank\NS(A)=1$.
\begin{center} 
   \def\bm{\bf}
   \small
   \begin{tabular}{rlll}
	 d & $\eps_{\rm lower}(d)$ & $\eps_{\rm upper}(d)$ & $\sqrt{2d}$ \\
	 \hline
                    \bm 1& 1.3228& \bm 1.3333& 1.4142 \\
		    \bm 2& 1.8708& --   & \bm 2 \\
                    3& 2.2912& 2.4000& 2.4494 \\
                    \bm 4& 2.6457& \bm 2.6666& 2.8284 \\
                    5& 2.9580& 3.1578& 3.1622 \\
                    6& 3.2403& 3.4285& 3.4641 \\
                    7& 3.5000& 3.7333& 3.7416 \\
		    \bm 8& 3.7416& --   & \bm 4 \\
                    9& 3.9686& 4.2352& 4.2426 \\
                    10& 4.1833& 4.4444& 4.4721 
   \end{tabular}
   \qquad\qquad
   \begin{tabular}{llll}
	 d & $\eps_{\rm lower}(d)$ & $\eps_{\rm upper}(d)$ & $\sqrt{2d}$ \\
	 \hline
                    11& 4.3874& 4.6903& 4.6904 \\
                    \bm 12& 4.5825& \bm 4.8000& 4.8989 \\
                    13& 4.7696& 5.0980& 5.0990 \\
                    14& 4.9497& 5.2913& 5.2915 \\
                    15& 5.1234& 5.4545& 5.4772 \\
                    16& 5.2915& 5.6470& 5.6568 \\
                    17& 5.4543& 5.8285& 5.8309 \\
		    \bm 18& 5.6124& --   & \bm 6 \\
                    19& 5.7662& 6.1621& 6.1644 \\
                    20& 5.9160& 6.3157& 6.3245
   \end{tabular}
\end{center}   
\end{remark}


\bigskip
\bigskip

   Thomas Bauer,
   Department of Mathematics,
   University of California, 
   Los Angeles, CA 90095--1555

   (E-mail: {\tt tbauer@math.ucla.edu})


\end{document}